# A Fast Lax-Hopf formula to solve the Lighthill-Whitham-Richards traffic flow model on networks

Michele D. Simoni[a],    Christian G. Claudel[b]

[a] Dept. of Civil, Architectural, and Environmental Engineering, University of Texas at Austin, m.simoni@utexas.edu

[b] Dept. of Civil, Architectural, and Environmental Engineering, University of Texas at Austin, christian.claudel@utexas.edu

**Abstract**

Efficient and exact algorithms are important for performing fast and accurate traffic network simulations with macroscopic traffic models. In this paper, we extend the semi-analytical Lax-Hopf algorithm in order to compute link inflows and outflows with the LWR model. Our proposed Fast Lax-Hopf algorithm has a very low computational complexity. We demonstrate that some of the original algorithm's operations (associated with the initial conditions) can be discarded, leading to a faster computation of boundary demand/supplies in network simulation problems, for general concave fundamental diagrams. Moreover, the computational cost can be further reduced for triangular Fundamental Diagrams and specific space-time discretizations. The resulting formulation has a performance comparable to the Link Transmission Model and, since it solves the original LWR model for a wide range of FD shapes, with any initial configuration, it is suitable to solve a broad range of traffic operations problems. As part of the analysis, we compare the performance of the proposed scheme to other well-known computational methods.

## 1. Introduction

Traffic flow models are commonly used to describe the propagation of traffic on transportation networks. Depending on the scale of the problem and on the type of traffic phenomena that need to be reproduced by the model, it is possible to identify three main classes of traffic flow models: microscopic, mesoscopic and macroscopic (Peeta and Ziliaskopoulos, 2001). In macroscopic models, traffic is modeled as a fluid stream described by a density and flow function, defined on all points of a road network, and for all times. Macroscopic models encode both the propagation of traffic on network links (resulting in macroscopic link models), as well as the splitting and merging of vehicle flows at junctions (resulting in junction models, or node models). One of the most commonly used macroscopic link model used in the literature is the the Lighthill–Whitham–Richards (LWR) model (Lighthill and Whitham, 1955; Richards, 1956). This model is based on two main assumptions: the conservation of vehicles and the existence of a univocal flow-density relationship (fundamental diagram). Assuming that links can be described by space-independent parameters (homogenous problem), the propagation of queues and shockwaves can be then modeled by means of a partial differential equation (PDE), known as the LWR PDE. The LWR model is often used for studies involving large simulations since it is relatively straightforward and robust, depending on a low number of model parameters that are easy to calibrate. Furthermore, its computational time that is independent of the number of vehicles to model (Wageningen-Kessel, 2016), unlike microscopic or mesoscopic models. Junction models have also been studied extensively to reproduce traffic behavior at merges/diverges (Daganzo, 1995), to investigate the propagation of kinematic waves (Garavello and Piccoli, 2006), and to identify general methods (Tampère et al., 2011; Flötteröd and Rohde, 2011; Jabari, 2016).

In the past two decades, a considerable number of numerical schemes have been proposed to solve the LWR model on networks, striving for higher computational efficiency and accuracy. The most popular ones include the Cell Transmission Model (CTM) (Daganzo, 1994), a particular case of Godunov discretization (Godunov, 1959), and the Link Transmission Model (LTM) (Yperman et al., 2006; Yperman, 2007), based on earlier work by Newell (Newell, 1993). Among the event-based numerical schemes, the Wave-front tracking methods (Bressan, 2000; Garavello and Piccoli,2006) reproduce the propagation of expansion waves and shocks using Riemann solvers and the Rankine Hugoniot formula (Baiti and Jenssen, 1998). Raadsen et al. (2016) propose another promising event-based algorithm suitable for large simulations, based on semi-analytical solutions of the LWR PDE. Event-based approaches can be very fast, however their efficiency and accuracy depend on the initial, boundary conditions and flux function of the problem to solve.

Alternative computational methods are based on the Hamilton-Jacobi formulation of the LWR model. Dynamic Programming (DP) methods (also referred as Variational Principle in the transportation literature) solve a network least cost problem (through DP) on space-time grid, resulting in the so-called Variational Method (Daganzo, 2005; Daganzo 2006). Alternatively, the Lax-Hopf (LH) algorithm (Lax, 1957; Hopf, 1969) uses a specific structure of the DP problem to decompose the solution into the minimum of a finite number of explicit functions, resulting in an exact computational method to compute the solution on a single link, and a faster computational time than DP. Unfortunately, the Lax-Hopf algorithm does not perform well over large time horizons and is in general slower than most of the traditional link models.

In this study, we propose a modification of the LH algorithm, referred as Fast Lax-Hopf (FLH) to solve the LWR model more efficiently while retaining the exactness of the LH. We show that its computational performance is comparable to the LTM, which is used, together with the CTM and the original LH, as a benchmark for our study.

Given known initial conditions in all links of the network, and given traffic demand and supply functions at the boundaries of the network, the objective of the present algorithm is to determine as quickly and precisely as possible the boundary conditions of each link. The boundary conditions are indeed a-priori unknown, and depend on the initial conditions across all links, the model used to describe the junctions, the model parameters and network layout. Once these boundary conditions are computed, the solution can be found at any point in space and time required by the particular problem (for example at a precise point of space and time where a measurement data point is generated for estimation problems, or at a given time horizon for forward simulation problems), by minimizing explicitly computed functions, using the classical LH algorithm. In the present article, we derive the FLH algorithm for general concave Fundamental Diagrams, that can be further simplified in the specific case of triangular Fundamental Diagram.

The FLH algorithm, which is particularly suitable for network simulations, requires lower number of operations than the original version of LH without compromising its accuracy. Furthermore, we show that the FLH can be further simplified for a Triangular Diagram, in specific space-time discretizations (CFL-like), while remaining exact. The resulting formulation shares similarities (in terms of the formulation and computational performance) with the well-known Link Transmission Model, though it is slightly slower than the latter. Nonetheless, the FLH algorithm can be used for any general concave Fundamental Diagram (FD), and arbitrary initial conditions, unlike the LTM.

The rest of this paper is organized as follows. First, we describe some of the main computational methods available, and discuss their advantages and drawbacks. We then derive the FLH algorithm using a set of theorems that simplify the original LH formulation. We also show that the FLH algorithm can be further simplified for triangular FDs and for some specific space-time discretizations, and show that the resulting algorithm has a similar (but not identical) expression to the original LTM. In the second part of the paper, we provide numerical validation of this algorithm by means of network traffic simulations, and comparisons with the original LH, CTM and LTM formulations. Finally, we present some considerations and conclusions based on the results.

## 2. Background: link models

Network simulation algorithms require a link model to reproduce traffic flow on each link. In this study, we focus on computational methods that solve the LWR model on each link of the domain (the parameters or even types of fundamental diagrams can change across the links of the network). In this section, after introducing the LWR model, we provide an overview of some of the main computational methods: LTM, CTM, Variational Theory, and LH (going through the details of each formulation is beyond the scope of the study).

### 2.1 The LWR model and the Hamilton-Jacobi PDE

For a given time $t$ and position $x$, we define the local traffic density $k(x,t)$ as the number of vehicles per unit length, and the instantaneous flow $Q(x,t)$ in vehicles per unit time. The conservation of vehicles on the highway is formulated as the following partial differential equation (PDE) (Lighthill and Whitman, 1956; Richards, 1956):

$$\frac{\partial k(t,x)}{\partial t} + \frac{\partial Q(t,x)}{\partial x} = 0 \tag{1}$$

In the Lighthill Whitham Richards (LWR) model, the Fundamental Diagram (FD) relates the flow and density; in this article, we consider general concave fundamental diagrams $Q(k)$ that are defined over some interval $[0, k_j]$, and for which $Q(0) = Q(k_j) = 0$. These general concave fundamental diagrams are the focus of Section 3.

An important particular case (studied in section 4) consists in the triangular FD (Daganzo, 1994). The FD is a positive and concave function defined on $[0, k_{jam}]$ where $k_{jam}$ is the maximal density (jam density). It ranges between $[0, q_{max}]$ where $q_{max}$ is the maximum flow (capacity). It is associated with derivatives $Q'(k) = v$ (free flow speed) for $k < k_c$ (critical density) and $Q'(k) = -w < 0$ (congested wave speed) for $k > k_c$. Hence, the triangular FD is defined as follows:

$$Q(k) = \begin{cases} v\,k & : \ 0 \leq k \leq k_c \\ -w\left(k - k_{jam}\right) & : \ k_c \leq k \leq k_{jam} \end{cases} \tag{2}$$

Since the triangular fundamental diagram is concave, it is continuous in the interior of its domain of definition, and therefore its parameters satisfy $vk_c = -w(k_c - k_{jam})$.

While the flow of traffic can be described by the density function $k(\cdot,\cdot)$, it can alternatively be described using the Moskovitz function $N(x,t)$ that expresses the cumulated vehicle count through a location $x$, at time $t$. The Moskowitz function (also called Cumulative number of vehicles function) is defined as follows. All vehicles on and entering the road link are labeled by increasing integers as they pass the entry point $x_0$ of a highway section, and are assumed not to pass each other. The Moskowitz function at location $x$ and time $t$ is defined as $N(x,t) = n$, where $n$ corresponds to the label of the vehicle closest to $x$ at time $t$. The derivatives of the Moskowitz function are related to the density and flow functions (Daganzo, 2006).

Replacing $k$ and $q$ with $N$ yields to the following Hamilton-Jacobi PDE (Newell, 1993; Daganzo, 2005a, 2006; Claudel and Bayen, 2010a, b):

$$\frac{\partial N(x,t)}{\partial t} - Q\left(-\frac{\partial N(x,t)}{\partial x}\right) = 0 \tag{3}$$

### 2.2 Computational methods

The LWR PDE is a first order hyperbolic scalar conservation law that can be solved using a number of computational methods.

In the CTM both time and space are discretized, as each link is divided into a given number of cells of size $\Delta x$. This size is constrained by the Courant-Friedrichs-Lewy (CFL) condition (Bretti et al., 2006), according to which, for a given time discretization $\Delta t$ the inequality $\Delta x \geq v\Delta t$ must hold, where $v$ is the free flow velocity. The CTM is essentially a Godunov discretization of the original LWR equation (when the flow-density relation is triangular or trapezoidal), and assumes that the density of vehicles in each cell is constant across space (Lebacque, 1995). For every time interval the number of vehicles leaving a given cell and entering in the cell immediately downstream is computed using the Godunov flux. The maximum number of vehicles that can fit into a cell is a function of the jam density. The CTM requires calculating flows for all the cells of the link in order to compute the upstream and downstream boundary conditions of this link. In addition, the CTM does not yield exact solutions to the LWR model in general, due to numerical diffusion errors (Leclercq et al., 2007). The discretization in cells leads to an approximation in the speed of shockwaves that can propagate over the network, and ultimately can yield considerable cumulated errors. Several extensions of the CTM have been proposed to model other properties of traffic, such as the capacity drop (Schreiter et al., 2010; Srivastava and Geroliminis, 2013), different shapes of the fundamental diagram (Lo, 1999), and to reduce the discretization error (Daganzo, 1999; Szeto, 2008). Although the CTM allows to fairly reproduce important traffic phenomena like the forming and propagation of queues, the spatial discretization of links represents a main limitation in terms of efficiency and accuracy (Gentile, 2010).

Instead, the LTM only requires time to be discretized. The main feature of this model based the simplified theory of Newell (1993a; b) consists in using the characteristic speeds (free-flow and congested flow) to derive the upstream and downstream boundary conditions. Recently, extensions of the original LTM formulation have been proposed to allow for larger time steps (Himpe et al., 2016) and to consider non-triangular FDs with capacity drops (Gun et al., 2017) and initial conditions (Gun, 2018). In recent years, the LTM has become very popular for the dynamic network loading (DNL) procedure within the dynamic traffic assignment (DTA), where simulations can involve thousands of links, and where the solution only needs to be computed on the link boundaries. However, a limitation of the LTM is that the solution cannot be computed inside each link, which makes it unsuitable to problems involving estimation and calculation of traffic indicators inside the links (e.g. in estimation, traffic optimization or control problems). In some specific situations, in which no expansion wave is present (for example in a constant initial density scenario), and for specific fundamental diagrams (triangular), the LTM allows computation of the solution inside the link, though this procedure does not converge towards the solution of the LWR model for general initial conditions or for general concave fundamental diagrams (see Section 5.3).

The Variational Theory introduced by Daganzo (2005) consists in applying Dynamic Programming to solve the Hamilton-Jacobi PDE (3) through the classical Lax-Hopf formula. The solution can equivalently be computed using the viability theory (Aubin et al., 2008). Both approaches are conceptually similar, with the exception that the viability approach allows more general (discontinuous) initial conditions to be considered, and allows the computation of lower-constrained solutions to the Hamilton-Jacobi PDE.

The Lax-Hopf algorithm exploits a particular structure of the Dynamic Programming problem used in the Variational Theory to compute the solutions more efficiently (and exactly) in the case where the fundamental diagram is space and time independent. In this situation, the solution can be obtained without discretizing the computational domain, and it corresponds to the minimum of a finite number of functions associated to the initial and boundary conditions. By definition, this method is analytical and yields exact results in simulations of single links. In the network simulations errors can occur due to the temporal discretization of the boundary conditions, since boundary conditions are not necessarily constant over a given time step.

Because it uses an additional structure of the DP problem, the Lax-Hopf algorithm is always faster than the Variational Theory (although the Variational Theory is more general since it cannot handle situations in which the fundamental diagram depends on space and/or time). Nevertheless, its computational performance is comparable to that of the CTM (Mazare et al., 2011) and thus offers no speed improvement over the abovementioned algorithms.

The FLH described in the following section allows one to compute solutions (at the boundaries) with lower computational requirements than the original LH. We achieve this by proving that some initial or boundary condition blocks appearing in the minimization problem (considered in the original LH) can be a-priori discarded, without affecting the results. Since these excluded blocks cannot theoretically influence the solution, the solution computed by this algorithm remains exact (for single link problems), as in Mazare et al. (2011). Once the sets of upstream and downstream boundary conditions has been derived through the FLH, they can be used to solve (Eq. 3) in any point of the computational domain, without relying on a computational grid, as in the original LH. This a particularly important aspect for estimation and control applications. For example, in estimation problems, one only needs to compute the solution on the space-time points corresponding to sensor measurements, which are in general considerably less than the total number of grid points (assuming an uniform grid in space and time). Similarly, in optimal control problems, the solution only needs to be computed on space-time points that are relevant for the computation of the objective function.

## 3. Fast Lax-Hopf Algorithm for computing solutions to the LWR model on networks

In this Section, we describe the main features of the LH algorithm used to compute the solutions of the LWR model semi analytically (Section 3.1). We then derive the FLH algorithm (Section 3.3) using a set of rules (Section 3.2) that can be used to reduce the number of calculations compared with the original LH algorithm.

### 3.1 The generalized Lax-Hopf Formula and boundary conditions

Let a value condition function $c(\cdot,\cdot)$ be defined. This value condition can encode for example initial and boundary conditions. Aubin et al. (2008) showed that the solution associated with the value condition $c(\cdot,\cdot)$, denoted here by $N_c(\cdot,\cdot)$, is the solution to the following optimization problem involving the value condition:

$$N_c = inf\{c(t - T, x - Tu) + TR(-u)\} \tag{4}$$
$$s.t.(u,T) \in [w,v] \times \mathbb{R}_+ \ and \ (t - T, x - Tu) \in Dom(c)$$

In the present article, the value condition $c(\cdot,\cdot)$ corresponds to initial, upstream and downstream boundary condition functions:

$$c(x,t) = \begin{cases} N_{ini}(x) & t = 0 \\ N_{up}(t) & x = x_0 \\ N_{down}(t) & x = x_n \end{cases} \tag{5}$$

The optimization problem (4) involves the function $R(\cdot)$, which is defined as the convex transform associated with the fundamental diagram $Q(\cdot)$:

$$R(u) = \sup_{k \epsilon [0,k_j]} (Q(k) - u \cdot k) \tag{6}$$

Equation (4) is well known in the Hamilton-Jacobi literature and often referred to as Lax-Hopf (LH) formula (Lax, 1973; Evans, 1998; Daganzo, 2006; Aubin et al., 2008; Claudel and Bayen, 2010 a,b). The convex transform $R(\cdot)$ is a convex, nonnegative and nonincreasing function.

The Lax-Hopf algorithm assumes that the initial and boundary conditions $c_{ini}(\cdot,\cdot), c_{up}(\cdot,\cdot)$ and $c_{down}(\cdot,\cdot)$ are piecewise linear (Mazare et al. 2011), and can thus be written as:

$$\begin{cases} c_{ini}(0,x) = c_{ini}^i(x) = -k_i x + b_i & if \ x_i \le x \le x_{i+1} \\ c_{up}(x_0,t) = c_{up}^j(t) = q_j t + d_j & if \ t_j \le t \le t_{j+1} \\ c_{down}(x_{n_{ini}},t) = c_{down}^j(t) = p_j t + c_j & if \ t_j \le t \le t_{j+1} \end{cases} \tag{7}$$

In this situation, the solutions associated with the $c_{ini}^i(x)$, $c_{up}^j(t)$ and $c_{down}^j(t)$ can be computed explicitly (Appendix I). The solution at any point $(t,x)$ of the space time domain can then be computed by taking the minimum of the solutions taken in $(t,x)$ and associated with each initial and boundary condition block. This comes from the inf-morphism property, initially derived in Aubin et al. (2008).

*3.2 Priority rules for computing the solution to the Hamilton-Jacobi equation*

The primary objective of the proposed algorithm is to quickly compute the outflows and inflows at every time step, by using a minimum number of operations, and maintaining exactness. Once the boundary conditions are known on all links, the solutions inside the computational domain can be found by minimizing a number of explicitly computed functions. The FLH algorithm speeds both the computation of the boundary conditions, and the computation of the solution inside the computational domain.

This algorithm relies on the specific structure of the partial solutions to the Hamilton-Jacobi PDE (Eq. 3). From (Claudel and Bayen, 2010a, b), the partial solutions associated with affine blocks are convex functions of $(t,x)$. Furthermore, (Daganzo 2005) showed that these solutions are Lipschitz continuous on their domain of definition for general concave fundamental diagrams.

In the present case, we consider a general mixed initial-boundary condition problem on a given stretch of highway limited by upstream and downstream boundaries. We also assume that the boundary conditions that apply on the domain are not known in advance, unlike in the LH case. These boundary conditions have to be computed at each time step through junction models relating the demands of the incoming links to the supplies of the outgoing links, across each junction. These junction models have the effect of coupling the solutions computed over adjacent links. Our objective is to compute the inputs to the junction models as fast as possible. These inputs are upstream demands and downstream supplies of each link (for a given time step). Once computed, these inputs can be used to determine the actual flows occurring at the upstream and downstream boundaries for the chosen time step, and use these results (initial and updated boundary conditions) to compute the solution at the subsequent time step. This allows one to iteratively compute the solution associated to (3) over transportation networks, for some given simulation time horizon T.

In the remainder of this article, we assume that the initial and boundary conditions are piecewise linear, resulting in piecewise constant density or flow blocks. More precisely, let the initial condition be expressed as a piecewise linear function, with each linear piece on intervals $(x_i, x_{i+1})$ defined by:

$$c_{ini}^i(x) = \begin{cases} -k_i x + b_i & : x_i \leq x \leq x_{i+1} \\ +\infty & : otherwise \end{cases} \tag{8}$$

where $i \in \{0, \dots, n_{ini} - 1\}$, with similar definitions for the upstream and boundary conditions (Eq. 7). As described in (Daganzo, 2006), the initial condition must satisfy some growth and continuity conditions:

$$0 \leq k_i \leq k_j \text{ for all } i \in \{0, \dots, n-1\} \tag{9}$$

$$-k_i x_i + b_i = -k_{i+1} x_i + b_{i+1}, \; \forall i \in \{1, \dots, n-1\} \tag{10}$$

Similar growth and continuity constraints apply for the upstream and boundary conditions, in particular the boundary flows are nonnegative and upper bounded by the link capacity.

By the inf-morphism property (Mazare et al, 2011), the solution $N(x,t)$ associated with the Hamilton-Jacobi PDE (3) can be computed at any point $(x,t)$ of the space-time domain using the following formula:

$$N(x,t) = \min\left( \min_{i,j,k} N_{c_{ini}^i}(x,t), N_{c_{up}^j}(x,t), N_{c_{down}^k}(x,t) \right) \tag{11}$$

To compute the downstream boundary block for a given time interval $[t, t + \Delta t]$, we first need to derive the demand of this particular link over the time interval $[t, t + \Delta t]$, defined by $d(t, t + \Delta t) = \frac{N\left(x_{n_{ini}}, t + \Delta t\right) - N(x_{n_{ini}}, t)}{\Delta t}$. The actual flow over the time interval $[t, t + \Delta t]$ is then determined using the other demand and of all links connected to this junction, through the chosen junction model.

Hence, assuming that $N(x_{n_{ini}}, t)$ is known, and using the classical LH algorithm (Mazare et al., 2011), we can compute $N(x_{n_{ini}}, t + \Delta t)$ as:

$$N\left(x_{n_{ini}}, t + \Delta t\right) = \min\left[ \min_{1 \leq j \leq n_{ini}} N_{c_{ini}^j}\left(x_{n_{ini}}, t + \Delta t\right), \min_{0 \leq k \leq \kappa} N_{c_{up}^k}\left(x_{n_{ini}}, t + \Delta t\right), N_{c_{down}^l}\left(x_{n_{ini}}, t + \Delta t\right) \right] \tag{12}$$

In (12), $\kappa$ is defined as $\kappa = \max_{i \in \mathbb{Z} \; s.t. \; x_0 + v_f \cdot (t + \Delta t - t_{i+1}) \geq x_{n_{ini}}} i$, and $l = \max_{i \in \mathbb{Z} \; s.t. \; t_{i+1} \leq t + \Delta t} i$.

For simplicity, we now assume that all boundary condition blocks are defined at regular time intervals (though the algorithm can be extended in a straightforward way for general time intervals), and thus, that $t_j = j \cdot \Delta t$, where $\Delta t$ is the time step considered. In this situation, we have that $\kappa = \lfloor \frac{t + \Delta t - \frac{x_{n_{ini}} - x_0}{v_f}}{\Delta t} \rfloor$ and $l = \lfloor \frac{t + \Delta t}{\Delta t} \rfloor$. In the original Lax-Hopf method, the process required to compute the downstream boundary condition block at time $t + \Delta t$ is shown in Figure 1.

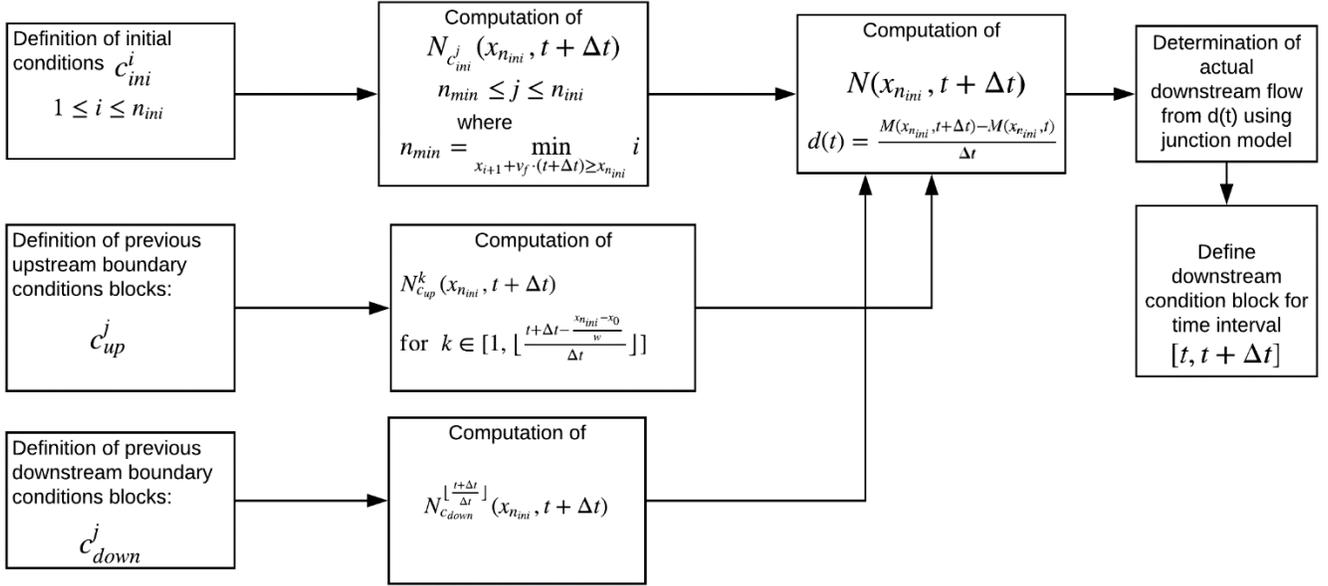

*Figure 1: required operations to determine the exiting flow (downstream) over the time interval $[t, +\Delta t]$ using the classical Lax-Hopf algorithm*

Equation (11) requires the minimization of $n_{ini} + \kappa + 1$ explicitly computed functions to derive the upstream supply of the link when $t + \Delta t \geq \frac{x_{n_{ini}} - x_0}{v_f}$. The objective of the Fast Lax Hopf algorithm is to decrease the required number of operations (in comparison to the Lax-Hopf algorithm), while still computing the average demand and supply functions exactly.

We now introduce a set of rules that allows one to reduce the number of required calculations with respect to the original LH algorithm.

**Theorem 1**: Let set of $n_{ini}$ initial conditions be defined as in (8). Let us further assume that $N_{c_{ini}^i}\left(x_{n_{ini}}, t'\right) \leq N_{c_{ini}^j}\left(x_{n_{ini}}, t'\right)$ for a time $t' \geq \frac{x_{n_{ini}} - x_{i+1}}{v_f}$, with $i < j$. Then:

$\forall\, t \geq t', N_{c_{ini}^i}\left(x_{n_{ini}}, t\right) \leq N_{c_{ini}^j}\left(x_{n_{ini}}, t\right)$

*Proof*: using the structure of the solution to initial conditions (Appendix I), we have that both $N_{c_{ini}^j}\left(x_{n_{ini}}, \cdot\right)$ and $N_{c_{ini}^i}\left(x_{n_{ini}}, \cdot\right)$ are defined on $[t', +\infty)$, continuous and convex functions. Hence, both functions have subderivatives (noted $\partial_-$), and are differentiable almost everywhere on their domain. These subderivatives can be computed from the expression of the initial solutions (Appendix I) as follows:

Let $v_i \in \partial_+ Q(k_i)$ (if $v_i \leq 0$ then only the third case of the below equation remains).

$$\partial_- N_{c_{ini}^i}\left(x_{n_{ini}}, t\right)$$

$$= \begin{cases} \{Q(k_i)\} & : \dfrac{x_{n_{ini}} - x_{i+1}}{v_i} \leq t \leq \dfrac{x_{n_{ini}} - x_i}{v_i} \\ \{R\left(\dfrac{x_{n_{ini}} - x_{i+1}}{t}\right)\} - \dfrac{x_{n_{ini}} - x_{i+1}}{t} \cdot \partial_- R\left(\dfrac{x_{n_{ini}} - x_{i+1}}{t}\right) : \dfrac{x_{n_{ini}} - x_{i+1}}{v_f} \leq t \leq \dfrac{x_{n_{ini}} - x_{i+1}}{v_i} \\ \{R\left(\dfrac{x_{n_{ini}} - x_i}{t}\right)\} - \dfrac{x_{n_{ini}} - x_i}{t} \cdot \partial_- R\left(\dfrac{x_{n_{ini}} - x_i}{t}\right) : t \geq \dfrac{x_{n_{ini}} - x_i}{v_i} \end{cases} \tag{13}$$

Using the Legendre-Fenchel inversion formula (Aubin Bayen Saint Pierre 2008), we have:

$$k \in \partial_- R\left(\frac{x_{n_{ini}} - x_{i+1}}{t}\right) \iff \frac{x_{n_{ini}} - x_{i+1}}{t} \in \partial_+ Q(k)$$

Hence, we have that $\left\{R\left(\frac{x_{n_{ini}} - x_{i+1}}{t}\right) - \frac{x_{n_{ini}} - x_{i+1}}{t} \cdot k\right\} = \{Q(k)\} \subset \{R\left(\frac{x_{n_{ini}} - x_{i+1}}{t}\right)\} - \frac{x_{n_{ini}} - x_{i+1}}{t} \cdot \partial_- R\left(\frac{x_{n_{ini}} - x_{i+1}}{t}\right)$ for $k \in \partial_- R\left(\frac{x_{n_{ini}} - x_{i+1}}{t}\right)$.

The same property can be applied to (13), allowing us to rewrite it as:

$$\partial_- N_{c_{ini}^i}\left(x_{n_{ini}}, t\right) = \begin{cases} \{Q(k_i)\} & : \dfrac{x_{n_{ini}} - x_{i+1}}{v_i} < t < \dfrac{x_{n_{ini}} - x_i}{v_i} \\ \left\{Q(k), k \in \partial_- R\left(\dfrac{x_{n_{ini}} - x_{i+1}}{t}\right)\right\} : \dfrac{x_{n_{ini}} - x_{i+1}}{v_f} < t < \dfrac{x_{n_{ini}} - x_{i+1}}{v_i} \\ \left\{Q(k), k \in \partial_- R\left(\dfrac{x_{n_{ini}} - x_i}{t}\right)\right\} : t > \dfrac{x_{n_{ini}} - x_i}{v_i} \end{cases} \tag{14}$$

We can rewrite $\partial_- N_{c_{ini}^i}\left(x_{n_{ini}}, t\right)$ as $\partial_- N_{c_{ini}^i}\left(x_{n_{ini}}, t\right) = Q(k_i(t))$ where $k_i(t)$ is the set-valued map defined by:

$$k_i(t) = \begin{cases} k_i = \partial_-(v_i) & : \dfrac{x_{n_{ini}} - x_{i+1}}{v_i} < t < \dfrac{x_{n_{ini}} - x_i}{v_i} \\ \partial_- R\left(\dfrac{x_{n_{ini}} - x_{i+1}}{t}\right) : \dfrac{x_{n_{ini}} - x_{i+1}}{v_f} < t < \dfrac{x_{n_{ini}} - x_{i+1}}{v_i} \\ \partial_- R\left(\dfrac{x_{n_{ini}} - x_i}{t}\right) : t > \dfrac{x_{n_{ini}} - x_i}{v_i} \end{cases} \tag{15}$$

It can be verified from this expression that $k_i(t) \leq k_j(t)$ when $i < j$, for $t > \frac{x_{n_{ini}} - x_{i+1}}{v_f}$. Indeed, since $R(\cdot)$ is convex, we have that $a \leq b \implies \partial_- R(a) \leq \partial_- R(b)$, in the sense of the interval order (partial order over the set of intervals of $\mathbb{R}$), and since $i < j$, we have $\frac{x_{n_{ini}} - x_{j+1}}{t} \leq \frac{x_{n_{ini}} - x_j}{t} \leq \frac{x_{n_{ini}} - x_{i+1}}{t} \leq \frac{x_{n_{ini}} - x_i}{t}$.

Hence, we have that $\partial_- N^i_{c_{ini}}(x_{n_{ini}}, t) \leq \partial_- N^j_{c_{ini}}(x_{n_{ini}}, t)$ since (in the sense of the interval order) for all times $t \geq t'$ if $i < j$.

Therefore, given that $N^i_{c_{ini}}(x_{n_{ini}}, t') \leq N^j_{c_{ini}}(x_{n_{ini}}, t')$, we have (by integration) that $N^i_{c_{ini}}(x_{n_{ini}}, t) \leq N^j_{c_{ini}}(x_{n_{ini}}, t)$.

**Theorem 2**: Let set of $n_{ini}$ initial conditions be defined as in (8). Let a set of upstream boundary conditions be defined as in Eq. 7. Let us assume that $N^i_{c_{ini}}(x_{n_{ini}}, t') \geq N^j_{c_{up}}(x_{n_{ini}}, t')$ for some $i \in [1, n_{ini}]$, for some time $t' > t_j + \frac{x_{n_{ini}} - x_0}{v_f}$. We have that $\forall t > t'$, $N^i_{c_{ini}}(x_{n_{ini}}, t') \geq N^j_{c_{up}}(x_{n_{ini}}, t')$.

*Proof*: given the structure of the solution to boundary conditions (Appendix I), and using a similar reasoning as in the proof or Theorem 1, we have that:

$$\partial_- N_{c_{up}j}(x_{n_{ini}}, t) = \begin{cases} Q(\rho_j) \ : \ x_0 + v_j(t - t_{j+1}) \leq x \leq x_0 + v_j(t - t_j) \\ \left\{ Q(k), k \in \partial_- R\left(\frac{x_{n_{ini}} - x_0}{t - t_j}\right) \right\} \ : \ x_0 + v_j(t - t_j) \leq x \leq x_0 + v_f(t - t_j) \\ \left\{ Q(k), k \in \partial_- R\left(\frac{x_{n_{ini}} - x_0}{t - t_{j+1}}\right) \right\} \ : \ x_0 \leq x \leq x_0 + v_j(t - t_{j+1}) \end{cases} \quad (16)$$

which can be rewritten as: $\partial_- N_{c_{up}j}(x_{n_{ini}}, t) = Q(k_j(t))$ where:

$$k_j(t) = \begin{cases} \rho_j \ : \ x_0 + v_j(t - t_{j+1}) \leq x \leq x_0 + v_j(t - t_j) \\ \partial_- R\left(\frac{x_{n_{ini}} - x_0}{t - t_j}\right) \ : \ x_0 + v_j(t - t_j) \leq x \leq x_0 + v_f(t - t_j) \\ \partial_- R\left(\frac{x_{n_{ini}} - x_0}{t - t_{j+1}}\right) \ : \ x_0 \leq x \leq x_0 + v_j(t - t_{j+1}) \end{cases} \quad (17)$$

It is straightforward to verify from both the above expression and Equation 15 in Theorem 1 that $k_j(t) \leq k_i(t)$ for all $t > t_j + \frac{x_{n_{ini}} - x_0}{v_f}$. This results again from the convexity of $R(\cdot)$, which implies that its subderivative is increasing.

Hence, we have that $\partial_- N^j_{c_{up}}(x_{n_{ini}}, t) \leq \partial_- N^i_{c_{ini}}(x_{n_{ini}}, t)$ for all times $t > t_j + \frac{x_{n_{ini}} - x_0}{v_f}$, and thus in particular for $t > t'$. Therefore, given that $N^i_{c_{ini}}(x_{n_{ini}}, t') \geq N^j_{c_{up}}(x_{n_{ini}}, t')$, we have that $\forall t > t'$, $N^i_{c_{ini}}(x_{n_{ini}}, t') \geq N^j_{c_{up}}(x_{n_{ini}}, t')$.

**Theorem 3**: Let a set of upstream boundary conditions be defined as in Eq. 7. Let us assume that $N_{c_{up}^j}\left(x_{n_{ini}}, t'\right) \leq N_{c_{up}^i}\left(x_{n_{ini}}, t'\right)$ for some, for some $i < j$, for some time $t' > t_j + \frac{x_{n_{ini}} - x_0}{v_f}$. We have that $\forall t > t'$, $N_{c_{up}^j}\left(x_{n_{ini}}, t'\right) \leq N_{c_{up}^i}\left(x_{n_{ini}}, t'\right)$.

*Proof*: The proof follows directly from the structure of the subderivative Eq. 16 and Eq. 17, remarking that if $i < j$, and $t > t_j + \frac{x_{n_{ini}} - x_0}{v_f}$, then $k_j(t) < k_i(t)$, and thus, $\partial_- N_{c_{up}^j}\left(x_{n_{ini}}, t\right) \leq \partial_- N_{c_{up}^i}\left(x_{n_{ini}}, t\right)$, which implies that $\forall t > t'$, $N_{c_{up}^j}\left(x_{n_{ini}}, t'\right) \leq N_{c_{up}^i}\left(x_{n_{ini}}, t'\right)$.

Equivalent properties can be derived for $x = x_0$ (upstream boundary), leading to the following theorems, which can be proved in a similar way as theorems 1, 2 and 3. For compactness we omit the proofs of these results, which are similar to the proofs outlined above.

**Theorem 4**: Let set of $n_{ini}$ initial conditions be defined as in (8). Let us further assume that $N_{c_{ini}^i}(x_0, t') \leq N_{c_{ini}^j}(x_0, t')$ for a time $t' \geq \frac{x_i - x_0}{w}$, with $j < i$. Then:

$\forall\, t \geq t', N_{c_{ini}^i}(x_0, t) \leq N_{c_{ini}^j}(x_0, t)$

**Theorem 5**: Let set of $n_{ini}$ initial conditions be defined as in (8). Let a set of downstream boundary conditions be defined as in Eq. 7. Let us assume that $N_{c_{ini}^i}(x_0, t') \geq N_{c_{down}^j}(x_0, t')$ for some $i \in [1, n_{ini}]$, for some time $t' > t_j + \frac{x_{n_{ini}} - x_0}{w}$. We have that $\forall t > t'$, $N_{c_{ini}^i}(x_0, t') \geq N_{c_{down}^j}(x_0, t')$.

**Theorem 6**: Let a set of downstream boundary conditions be defined as in Eq. 7. Let us assume that $N_{c_{up}^j}(x_0, t') \leq N_{c_{up}^i}(x_0, t')$ for some, for some $i < j$, for some time $t' > t_j + \frac{x_{n_{ini}} - x_0}{w}$. We have that $\forall t > t'$, $N_{c_{up}^j}(x_0, t') \leq N_{c_{up}^i}(x_0, t')$.

### 3.3 The Fast Lax-Hopf Algorithm

The Fast Lax-Hopf algorithm leverages the Theorems outlined in section 3.2 to reduce the number of calculations required to determine the solutions at a given boundary (upstream or downstream) of the domain. The Fast Lax-Hopf algorithm is summarized in Figure 2 below, for the computation of the downstream demand over time, over a single road link.

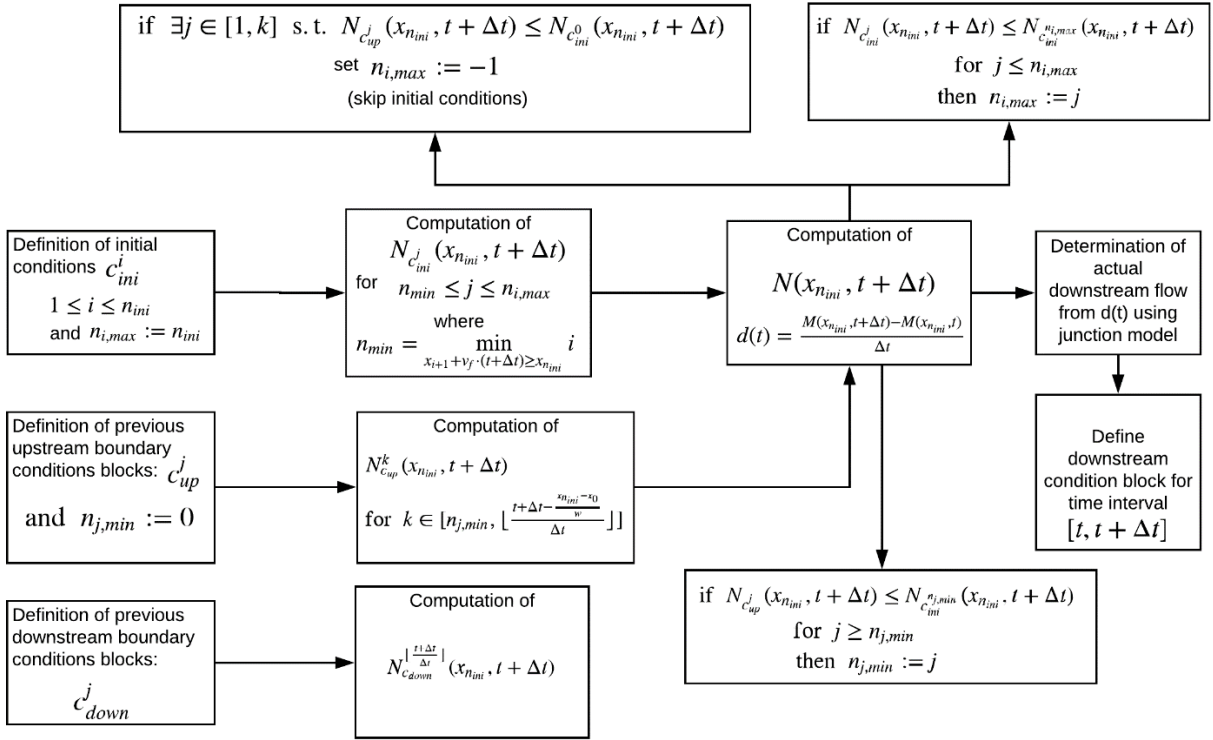

*Figure 2: required operations to determine the exiting flow (downstream) over the time interval [t,+Δt] using the Fast Lax-Hopf algorithm*

Given its similar structure as the Lax-Hopf algorithm (only with less operations), the FLH algorithm has a complexity that is upper bounded by that of the Lax-Hopf algorithm, while still retaining the exactness of the former for single link problems with piecewise constant demand and supply functions. Note that exactness is lost (as with all other algorithms available, except possibly wave-front tracking) on network problems since boundary demand and supply functions are not piecewise constant in general.

## 4.   Fast Lax-Hopf Algorithm for triangular fundamental diagrams

In this section, the Fast Lax-Hopf Algorithm formulation for triangular fundamental diagrams (FDs) is presented (Section 4.1). A more specific formulation is derived for particular situations where initial conditions have constant size and the time step satisfies the CFL condition (Section 4.2). Its computational complexity in comparison with other known numerical approaches is discussed in Section 4.3

*4.1 Specific formulation of triangular FDs*

In the specific case of a triangular fundamental diagram, the convex transform $R(\cdot)$ is affine:

$$\forall u \in [-w, v], \qquad R(u) = k_c(v - u) \tag{18}$$

and the solutions to affine initial, upstream or downstream boundary conditions are piecewise linear, as shown in Mazare et al. (2011), and can be written explicitly. In particular, the solution at any arbitrary point $(t, x)$ depends only upon at one specific (predictable) upstream boundary condition block, and one downstream boundary condition block, leading to the following:

$$N(x,t) = \min(min_{i \leq j \leq n_{ini}} N_{c_{ini}^j}(x,t), N_{c_{down}^{\left\lfloor \frac{t}{\Delta t} - \frac{x_{n_{ini}} - x_0}{w \Delta t} \right\rfloor}}(x,t), N_{c_{up}^{\left\lfloor \frac{t}{\Delta t} - \frac{x_{n_{ini}} - x_0}{v_f \Delta t} \right\rfloor}}(x,t)) \tag{19}$$

Furthermore, the number of required operations required to compute the solution at an arbitrary point $(t, x)$ of the computational domain can also be reduced as follows:

**Corollary 1**: Let a set of $n_{ini}$ initial conditions be defined as in (8), with Lipschitz continuity constraints (9) and (10). Let us further assume that $N_{c_{ini}^j}(x, t_s) \leq N_{c_{ini}^i}(x, t_s)$ for a time $t_s \geq \frac{x_{i+1} - x}{w}$, with $i < j$. Then:

$$\forall t \geq s, N_{c_{ini}^j}(x,t) \leq N_{c_{ini}^i}(x,t) \tag{20}$$

Proof: using the structure of the solutions $N_{c_{ini}^j}(x, t)$, we have that $N_{c_{ini}^j}(x, t) \leq N_{c_{ini}^j}(x, t_s) + (t_s - t)v\, k_c$ if $t_s \geq \frac{x_{i+1} - x}{w}$ and $i < j$, irrespective of the value of $k_j$. We also have that $N_{c_{ini}^i}(x, t) = N_{c_{ini}^i}(x, t_s) + (t_s - t)v\, k_c$. Since $N_{c_{ini}^i}(x, t_s) \leq N_{c_{ini}^j}(x, t_s)$, we have that $\forall t \geq t_s, N_{c_{ini}^j}(x, t) \leq N_{c_{ini}^i}(x, t)$.

This theorem implies that inside the computational domain, if the solution associated to a particular initial condition piece $j$ is lower than the solution associated with another initial condition piece $i$ (with $i < j$), for a location $x$ and time $t_s$ such that $t_s \geq \frac{x_{i+1} - x}{w}$, then the solution associated with piece $i$ cannot influence the solution (at the same location) at subsequent times.

**Corollary 2**: Let a set of $n_{ini}$ initial conditions be defined as in (8), with Lipschitz continuity constraints (9) and (10). Let us further assume that $N_{c_{ini}^j}(x, t_V) \leq N_{c_{ini}^i}(x, t_V)$ for some $t_V \geq \frac{x - x_i}{v_f}$, with $i > j$. Then:

$$\forall\, t \geq t_V, N_{c_{ini}^j}(x_0, t) \leq N_{c_{ini}^i}(x_0, t) \tag{21}$$

Proof: using the structure of the solutions $N_{c_{ini}^j}(x,t)$ , we have that $N_{c_{ini}^i}(x,t) = N_{c_{ini}^j}(x,t_V) + (t_V - t)v\, k_c$ if $t_V \geq \frac{x-x_i}{v}$ and $i > j$, irrespective of the value of $k_j$. We also have that $N_{c_{ini}^j}(x,t) \leq N_{c_{ini}^j}(x,t_V) + (t_V - t)v\, k_c$. Hence, we have that $\forall\, t \geq t_V, N_{c_{ini}^j}(x,t) \leq N_{c_{ini}^i}(x,t)$.

This result similarly allows us to exclude a priori some terms from (19), and can be used to speed up computations inside the computational domain. Hence, with the above rules, the computation of the solution at any point of the computational domain can be further simplified as:

$$N(x,t) = \min(\min_{j \in S(x,t)} N_{c_{ini}^j}(x,t), N_{c_{down}}^{\lfloor \frac{t}{\Delta t} - \frac{x_{n_{ini}} - x_0}{w\,\Delta t} \rfloor}(x,t), N_{c_{up}}^{\lfloor \frac{t}{\Delta t} - \frac{x_{n_{ini}} - x_0}{v_f\,\Delta t} \rfloor}(x,t)) \tag{22}$$

where $S(x,t)$ is the set of initial conditions indices used for the computation $N(x,t)$, which is updated using Corollary 1 and Corollary 2.

Note that $S(x,t)$ depends on the structure of the initial conditions, and is difficult to compute a-priori, though it can be iteratively computed on a computer using Corollary 1 and Corollary 2. If, in addition, the solution is computed at the boundaries of the domain, and the discretization of the initial conditions follows a CFL-type condition, $S(x,t)$ can be computed straightforwardly.

*4.2 Formulation for specific spatio-temporal discretizations*

In this section, we further assume that the domains of the initial condition satisfy $x_i = x_0 + i\Delta x$ (where $i \in N$), that is, that the initial conditions are piecewise constant on domains of constant size $\Delta x$. We also assume that the space and time steps satisfy a CFL-type condition: $\Delta t \leq \frac{\Delta x}{v}$ . In this situation, we can prove the two following results (for the upstream boundary, the downstream boundary case being similar), which further simplify the computation of the solution at the upstream and downstream boundaries:

**Corollary 3**: Let a set of $n_{ini}$ initial conditions be defined as in (8), with Lipschitz continuity constraints (9) and (10). Let us further assume that $x_i = x_0 + i\Delta x$ and $\Delta t \leq \frac{\Delta x}{v}$. For any discrete time $t = i \cdot \Delta t,\ i \in N$ we have that:

$$N(x_0,t) = \begin{cases} \min\left(N_{c_{ini}^l}(x_0,t), N_{c_{ini}^{l-1}}(x_0,t), N_{c_{up}}{}^j(x_0,(i-1)\Delta t) + v \cdot k_c \cdot \Delta t\right) \; if \; t \leq \frac{x_{n_{ini}}-x_0}{w} \\ \min\left(N_{c_{up}}{}^j(x_0,(i-1)\Delta t) + v \cdot k_c \cdot \Delta t, N_{down}{}^k(x_0,t)\right) \qquad\qquad else \end{cases}$$

(23)

where $j = i - 1$ , $k = \left\lfloor \frac{t-\frac{x_{n_{ini}}-x_0}{w}}{\Delta t} \right\rfloor$, $l = \left\lfloor \frac{wt}{\Delta x} \right\rfloor$

Proof: The first case corresponds to the situation where only initial components and upstream boundary condition components can influence the upstream condition ($t \leq \frac{x_{n_{ini}}-x_0}{w}$). In this situation, we have that $N_{c_{ini}^k}(x_0,t) = +\infty$ if $k > l$. Hence, we can write that $N(x_0,t) = \min\left(N_{c_{ini}^0}(x_0,t), ..., N_{c_{ini}^{l-1}}(x_0,t), N_{c_{ini}^l}(x_0,t), N_{c_{up}}{}^j(x_0,(i-1)\Delta t) + v \cdot k_c \cdot \Delta t\right)$. However, by the structure of the initial condition solution components (12), we have that for any $k \in [0, l-2]$, $N_{c_{ini}^k}(x_0,t) = N_{c_{ini}^k}(x_0,(i-1)\Delta t) + vk_c(t - (i-1)\Delta t)$. By the inf-morphism property $N_{c_{ini}^k}(x_0,(i-1)\Delta t) \geq N(x_0,(i-1)\Delta t)$, and thus, since $N(x_0,t) \leq N(x_0,(i-1)\Delta t) + k_c v(t - (i-1)\Delta t)$, we have that $N(x_0,t) \leq N_{c_{ini}^k}(x_0,t)$ for any $k \in \{0, ..., l-2\}$, which shows that only $N_{c_{ini}^{l-1}}$ or $N_{c_{ini}^l}$ can influence the solution in $(x_0,t)$. The proof of the second case is similar.

**Corollary 4**: Let a set of $n_{ini}$ initial conditions be defined as in (8), with Lipschitz continuity constraints (9) and (10). Let us further assume that $x_i = x_0 + i\Delta x$ and $\Delta t \leq \frac{\Delta x}{v}$. For any discrete time $t = i \cdot \Delta t$ ($i \in N$) we have that:

$$N(x_{n_{ini}},t) = \begin{cases} \min\left(N_{c_{ini}^l}(x_{n_{ini}},t), N_{c_{ini}^{l+1}}(x_{n_{ini}},t), N_{c_{down}}{}^j(x_{n_{ini}},(i-1)\Delta t) + v \cdot k_c \cdot \Delta t\right) \; if \; t \leq \frac{x_{n_{ini}}-x_0}{v} \\ \min\left(N_{c_{down}}{}^j(x_{n_{ini}},(i-1)\Delta t) + v \cdot k_c \cdot \Delta t, N_{up}{}^k(x_{n_{ini}},t)\right) \qquad\qquad else \end{cases}$$

(24)

where $j = i - 1$ , $k = \left\lfloor \frac{t-\frac{x_{n_{ini}}-x_0}{v}}{\Delta t} \right\rfloor$, $l = \left\lfloor \frac{vt}{\Delta x} \right\rfloor$

Proof: The first case corresponds to the situation where only initial components and upstream boundary condition components can influence the upstream condition ($t \leq \frac{x_{n_{ini}}-x_0}{v}$). In this situation, we have that $N_{c_{ini}^k}(x_{n_{ini}},t) = +\infty$ if $k > l$. Hence, we can write that $N(x_{n_{ini}},t) = \min\left(N_{c_{ini}^0}(x_{n_{ini}},t), ..., N_{c_{ini}^{l-1}}(x_{n_{ini}},t), N_{c_{ini}^l}(x_{n_{ini}},t), N_{c_{up}}{}^j(x_{n_{ini}},(i-1)\Delta t) + v \cdot k_c \cdot \Delta t\right)$. However, by the structure of the initial condition solution components (12), we have that for any $k \in [0, l-2]$, $N_{c_{ini}^k}(x_{n_{ini}}, t) = N_{c_{ini}^k}(x_{n_{ini}}, (i-1)\Delta t) + vk_c(t - (i-1)\Delta t)$. By the inf-morphism property $N_{c_{ini}^k}(x_{n_{ini}}, (i-1)\Delta t) \geq N(x_{n_{ini}}, (i-1)\Delta t)$, and thus, since $N(x_{n_{ini}}, t) \leq N(x_{n_{ini}}, (i-1)\Delta t) + k_c v(t - (i-1)\Delta t)$, we have that $N(x_{n_{ini}}, t) \leq N_{c_{ini}^k}(x_{n_{ini}}, t)$ for any $k \in \{0, \dots, l-2\}$, which shows that only $N_{c_{ini}^{l+1}}$ or $N_{c_{ini}^l}$ can influence the solution in $(x_{n_{ini}}, t)$. The proof of the second case is similar.

The above results (and their downstream boundary condition counterparts) imply that, when computing the upstream and downstream conditions in the initial phase of the computation, the associated solutions can be computed on just two consecutive blocks (Figure 3a-b). Furthermore subsequent computations of the solutions at the upstream and downstream boundaries (outside of the area of influence of the initial conditions) can be reduced to those in the classical LTM formulation. The FLH scheme thus computes the boundary conditions with a slightly higher computational cost as the LTM during the initial phase of the simulation (in the domain of influence of the initial condition), requiring two operations to account for the initial conditions instead of one. For subsequent times both formulations (LTM and FLH) are identical, and thus have the same computational cost.

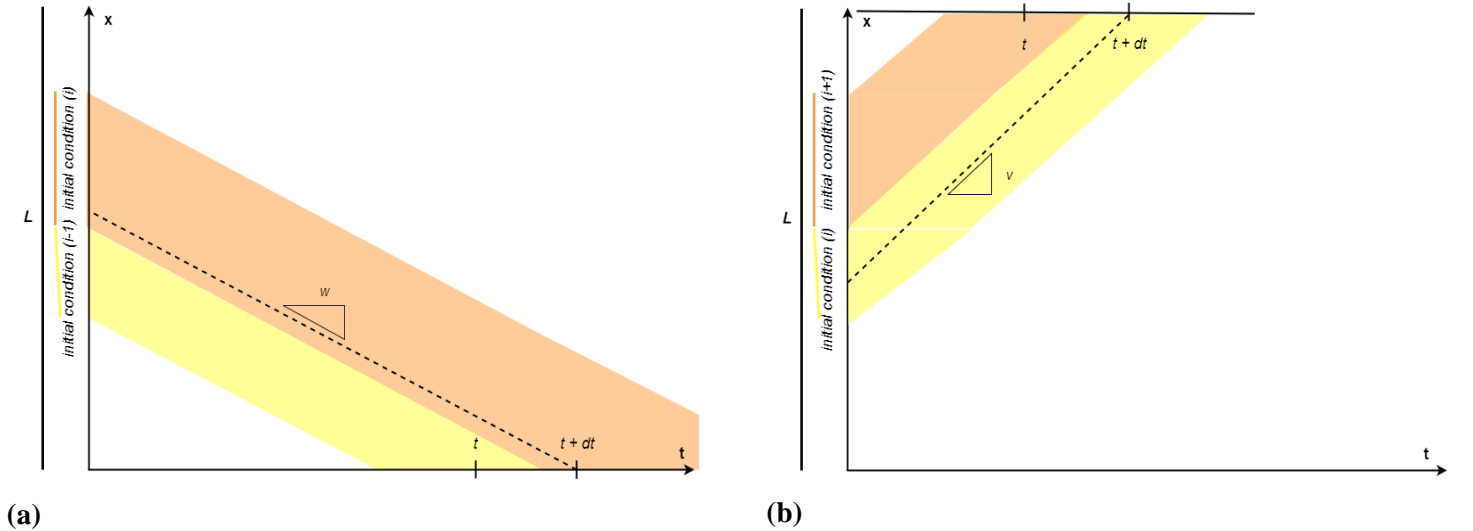

**(a)**             **(b)**

*Figure 3: Initial conditions considered for computation of flows upstream (a) and downstream (b) according to Theorem 5 and Theorem 6*

### 4.3 Comparison of the computational complexity of the FLH with other algorithms

The computational time required by the different algorithms outlined in the previous section depend on the type of problem that needs to be solved. In this section, we once again assume a triangular FD for simplicity (a non-triangular diagram would lead to non-convergence for the LTM algorithm in general). We consider two different problems:

1. Computing the solution to the LWR model at the **boundaries** of each link within a road network, with temporal step $\Delta t$, over some time horizon $T$.

2. Computing the solution to the LWR model inside the link, on a **uniform grid** of spatial resolution $\Delta x$ and temporal resolution $\Delta t$, over some time horizon $T$.

The first problem is typically encountered in forward simulations involving network loading, or network control when the objective function depends only upon the state computed at the boundaries of each computational domain. In contrast, the second problem is common in the applications such as traffic estimation (Cristiani et al., 2011), control (Ferrara et al., 2015), and estimation (Work et al., 2008). In some of these problems the solution only needs to be computed at specific points of the space-time grid, which are known in advance (for example, in estimation problems the solution only needs to be computed on points corresponding to sensors' locations). In this type of situation, the FLH can be used to solve the solution just at these specific locations without relying on the entire grid like in the second problem, unlike CTM and Variational Theory.

The computational performance of all algorithms is illustrated in Table 1 below, where $n_{ini}$ represents the number of initial conditions (or the number of grid points in the $x$ axis), and $n_t$ represents the number of time steps ($n_t = \frac{T}{\Delta t}$).

In the first problem, although the algorithms differ by their accuracy and computational cost, they all converge to the true solution to the LWR model when $\Delta t \to 0$ (and when both $\Delta t \to 0$ and $\Delta x \to 0$ for the CTM and dynamic programming). The LTM is the fastest algorithm, requiring 2 calculations per time step. In contrast, the CTM requires at least $4 \cdot n_{ini}$ calculations (computing demand and supplies, and computing the Godunov flux) per time step. The LH algorithm and DP both require on the order of $n_{ini}$ computations per time step (less during the first time steps), and are thus not significantly improving over the CTM. In contrast, the FLH algorithm requires 3 calculation per time step when $t \leq \frac{x_{n_{ini}} - x_0}{w}$, and 2 calculations per time step when $t > \frac{x_{n_{ini}} - x_0}{w}$. It thus has a computational complexity comparable to that of the LTM.

The second problem can be solved by all algorithms except the LTM, which is designed by definition for computations only at links' boundaries and it is not convergent whenever the initial condition contains expansion waves (see Section 5.3 for further explanations). Similarly to the first problem, the CTM requires $4 \cdot n_{ini}$ calculations per time step. DP methods require a number of calculations on the order of $n_{ini}^2$ per time step, while the classical LH algorithm requires on the order of $n_{ini} \cdot (n_{ini} + 2)$ calculations per step, which is similar to the DP. In contrast, the FLH algorithm requires less than $3 \cdot n_{ini}$ calculations per time step, which is a considerable improvement, and on par with the classical CTM. Note that in practice the FLH can be considerably faster than the CTM when the solution does not have to be computed in all cells (for example, in most of estimation, control or optimization problems).

*Table 1:* computational performance and accuracy of different algorithms (Triangular Fundamental Diagram)

| Numerical scheme | DP | LH | FLH | CTM | LTM |
|---|---|---|---|---|---|
| **Computational complexity (1)** | $\sim n_{ini} \cdot n_t$ | $\sim n_{ini} \cdot n_t$ | $2 \cdot n_t \leq$ $\leq 3 \cdot n_t$ | $\sim n_{ini} \cdot n_t$ | $\sim 2 \cdot n_t$ |
| **Computational complexity (2)** | $\sim n_{ini}^2 \cdot n_t$ | $\sim n_{ini}^2 \cdot n_t$ | $\lesssim n_{ini} \cdot n_t$ | $\sim n_{ini} \cdot n_t$ | Not convergent |
| **Accuracy** | Convergent | Exact on single links problems | Exact on single links problems | Convergent | Convergent on boundary conditions simulation problems |

## 5. Numerical Implementation

In this section we numerically implement the FLH algorithm presented in Section 3, and compare it to the LH, LTM and CTM, when possible. We show that the FLH algorithm has favorable characteristics in comparison with existing algorithms, particularly when solutions do not need to be calculated everywhere. In Section 5.1, we describe the computation of an initial condition problem over a single link, with a Greenshields FD. In Section 5.2, we focus on the specific case of a triangular FD, in a network simulation problem.

*5.1 Single link case*

The computation of solutions over a single link, with free upstream and downstream boundary conditions (for simplicity), is illustrated here. We consider a Greenshields FD given by $Q(k) = \frac{v_f}{k_j} k \cdot (k_j - k)$ where $v_f = 1$ and $k_{jam} = 4$ (arbitrary units), and a piecewise constant initial condition defined over regular intervals as follows:

$$k_i(x) = \begin{cases} 1.9 \ if \ x \in [0,40] \\ 3.0 \ if \ x \in [40,80] \\ 0.1 \ if \ x \in [80,120] \\ 3.7 \ if \ x \in [120,160] \\ 2.6 \ if \ x \in [160,200] \\ 4.0 \ if \ x \in [200,240] \\ 3.3 \ if \ x \in [240,280] \\ 0.4 \ if \ x \in [280,320] \\ 1.0 \ if \ x \in [320,360] \\ 0.3 \ if \ x \in [360,400] \end{cases}$$

The solution associated with a time horizon of 400 time steps is computed, and the corresponding density function is shown in Figure 4.

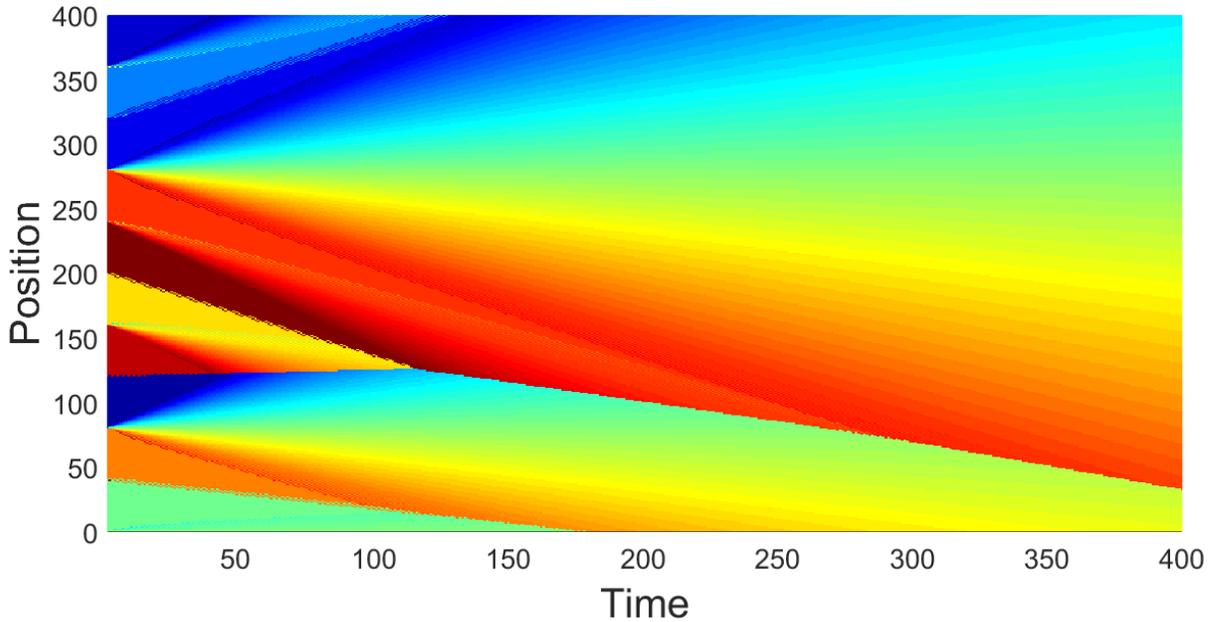

*Figure 4: Space-time-density diagram associated with the given initial conditions*

We now compute the domains of influence associated with this solution. The domains of influence $K_i, i \in [1,10]$ corresponds to the subset of computational domain defined by $K_i = \{(x,t) \in [x_0, x_{n_{ini}}] \times \mathbb{R}_+ \mid N(x,t) = N_{c_{ini}^i}(x,t)\}$. Note that these sets do not form a partition of $[x_0, x_{n_{ini}}] \times \mathbb{R}_+$ since some areas of the computational domain are minimal for multiple different initial condition blocks at the same time. Figure 5 provides an illustration of Theorem 1 and Theorem 4 for the computation of the downstream demand and upstream supply respectively. A consequence of Theorem 1 is that the indices of the initial condition blocks minimizing $N$ are decreasing functions of time on the downstream boundary, and increasing functions of time on the upstream boundary (Theorem 4), which can be clearly seen in this figure. Using this rule, it is possible to significantly reduce the computational time required to determine the solution at the boundaries: for example, past time t=50, the block 10 can be discarded from the computation of the downstream demand, and past t=150, the blocks 8, 9 and 10 can be discarded from the computation of the downstream demand. The same applies for block 1 and the computation of the upstream supply according to Theorem 4.

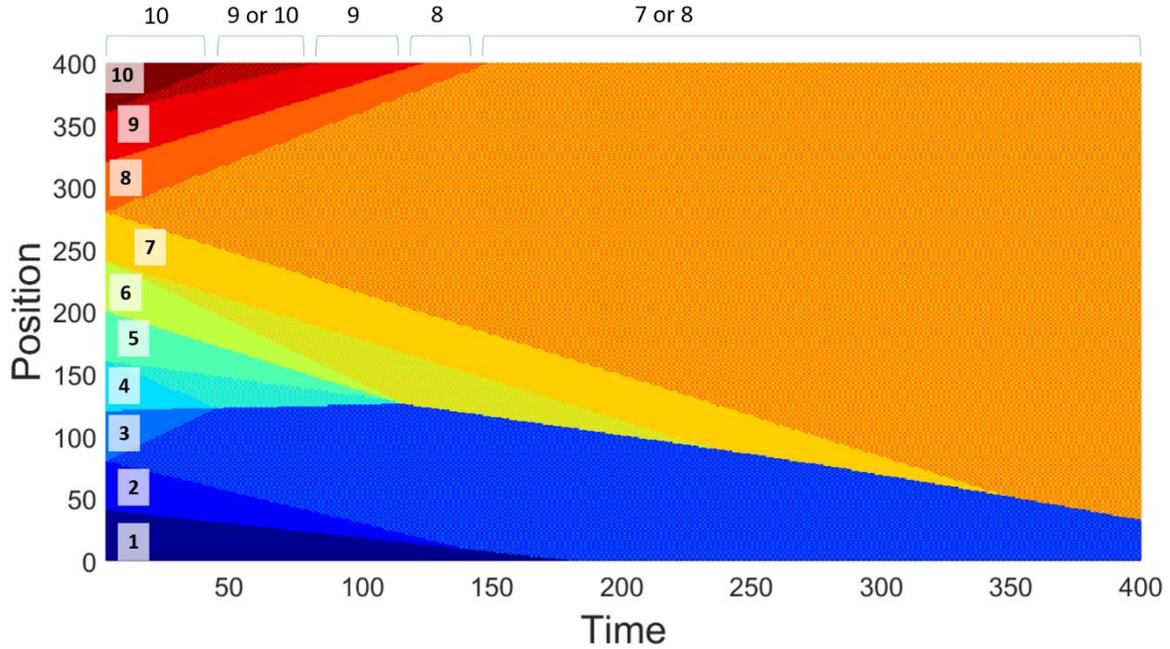

*Figure 5: Corresponding domains of influence $K_i$ associated with the given initial conditions. Note that the domains of influence can sometimes have a nonempty intersection, since multiple initial condition blocks may be minimal at the same time.*

### 5.2 Network case study

In the case of network simulations, the higher accuracy and speed of methods like the FLH and LTM compared to the CTM becomes apparent. As an example we show the simulation results of a 5 link highway network (Figure 6) composed of a 3-lane major highway section (Links 1,2 and 3), a 2-lane off-ramp (Link 4) and a 2-lane on-ramp (Link 5). A triangular fundamental diagram with capacity $q_{max}$=0.556 veh/s/lane, free-flow speed $v$=30 m/s and jam density of $k_{jam}$=0.1297 veh/m/lane is adopted for all three models. All links are characterized by initial free-flow density $k_1$=0.004 veh/m, with the exception of Link 2, that is characterized by two initial condition blocks associated with densities $k_1$=0.004 veh/m in the downstream half and $k_2 = 0.01$ veh/m in the upstream half of the link.

In order to model traffic throughout intersections, there is need of a generic macroscopic node model that respects some critical conditions: satisfaction of links' capacity constraints; conservation of flows; satisfaction of demand distribution constraints; maximization of flows (vehicles should proceed if there is available supply downstream); satisfaction of invariance principle (if the flows are restricted by demands, solutions cannot vary by increasing supplies and vice versa); and non-simultaneity of conflicting flows. In this study, we adopt the "I-HFS algorithm" by Jabari (2016), which respects the abovementioned properties and efficiently derives solutions by staging movements according to any arbitrary priory rules.

As can be seen in Figure 7, the results of the LTM and LH algorithms are close to the exact solution to the problem, while the solution computed by the CTM exhibits significant errors. We then compare the performance of all three schemes (CTM, LTM and FLH) in computing the solutions at the boundaries of each link of the network, averaged over random initial conditions, and random boundary demand and supplies at the edge of the road network. The results are averaged over 100 simulations, where the initial condition densities, demand and supply flows are drawn independently from uniform distributions. The average Root Mean Square error (RMSE) of all three schemes is shown in Figure 8. As can be seen from this Figure, FLH outperforms both the CTM and LTM in terms of error.

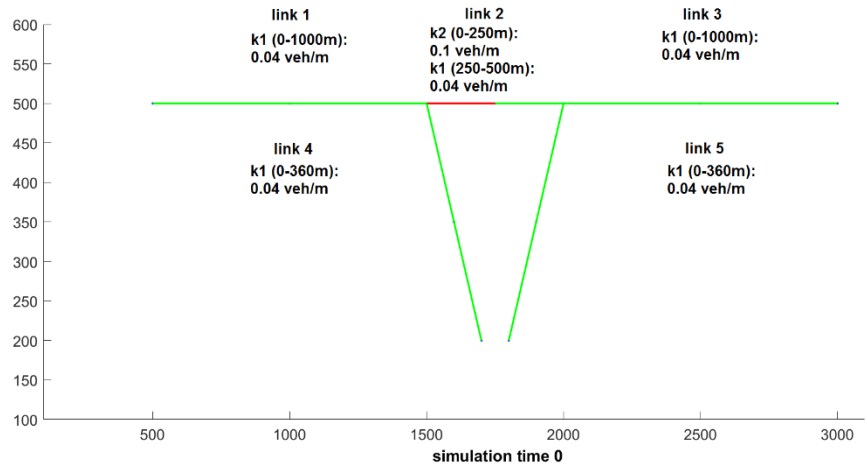

Figure 6: Simulation of the highway network at t=0 seconds.

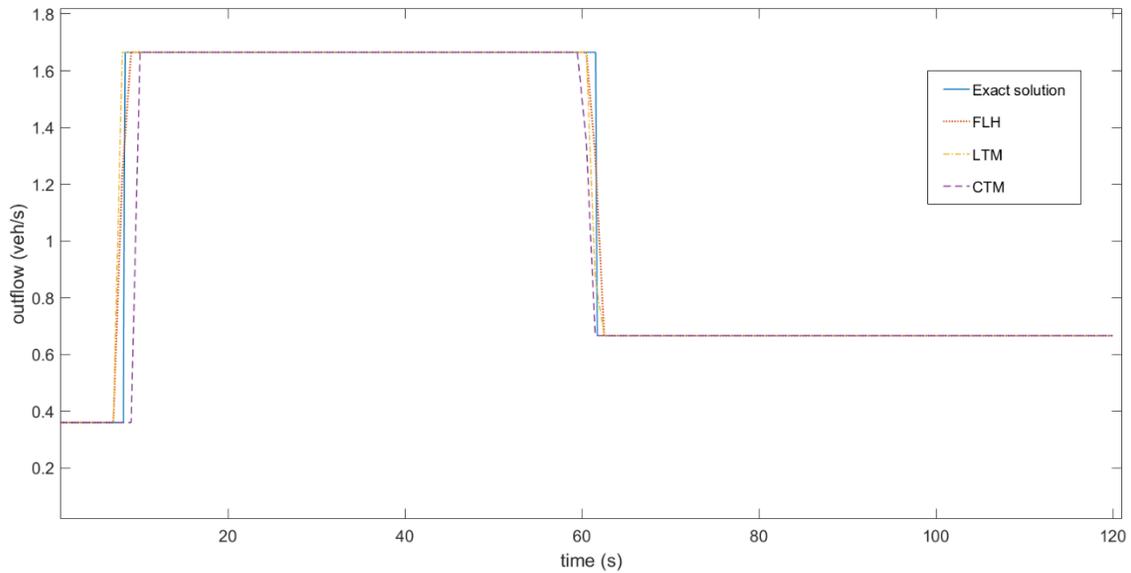

Figure 7: Comparison of the outflows of Link 2 with the three methods (using a time step of 1 second)

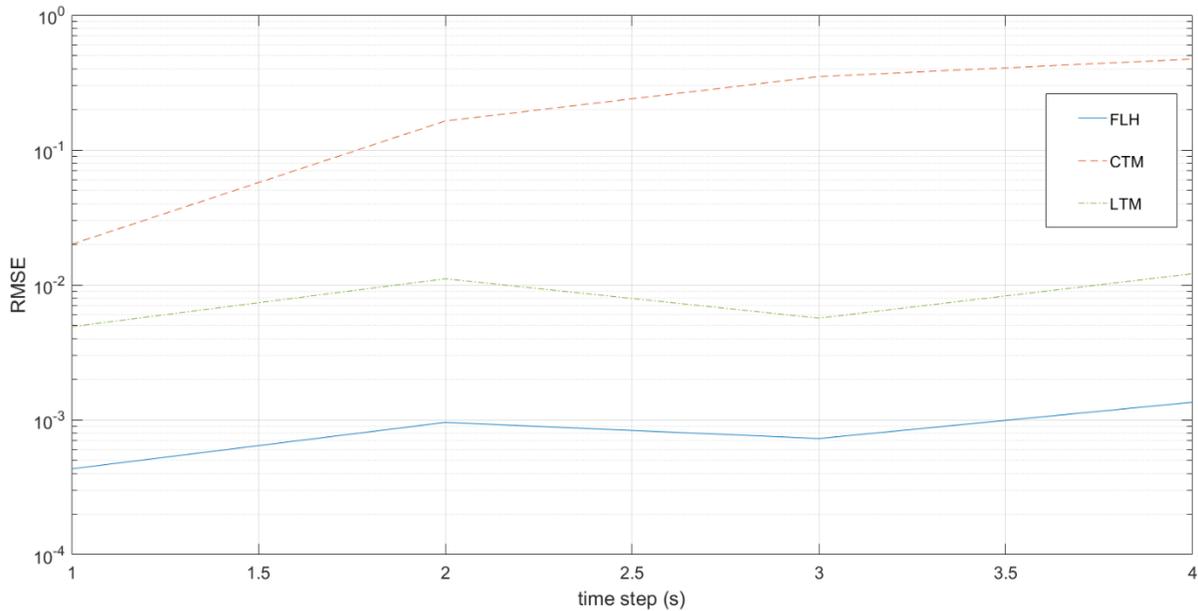

Figure 8: Accuracy of outflows calculated with the three different methods according to increasing time step

The favourable computational time properties of the FLH appear also in large network simulations. In order to demonstrate the scalability of the algorithm, we present and discuss the results of its application to a subset of the Austin downtown network (Figure 9). The network is characterized by 201 links and 110 nodes. Streets have between 2 to 3 lanes and the majority of the intersections is signalized (about 90%). For simplicity, in this study we only model green/red phases and we adopt the same triangular fundamental diagram for all links with: $q_{max}$=0.4625 veh/s, $v$=12.5 m/s, and $k_{jam}$=0.1295 veh/m.

We report in Table 2 and Table 3, the (average) computation times obtained for increasing simulation horizons, using different time steps for the four models. The simulations were performed on Matlab running on a laptop with a 2.8 GHz processor. The results are consistent with those obtained for single-link simulations. FLH and LTM have comparable performances when the initial conditions are explicitly considered in the computation. The CTM, for larger time steps (e.g. 5 seconds) is equivalent to the other methods since the links of this network are relatively short (resulting in a low numbers of cells). The numerical approximation of CTM, however, amplifies on large networks, leading to significant divergence from the exact solution after relatively short simulation horizons.

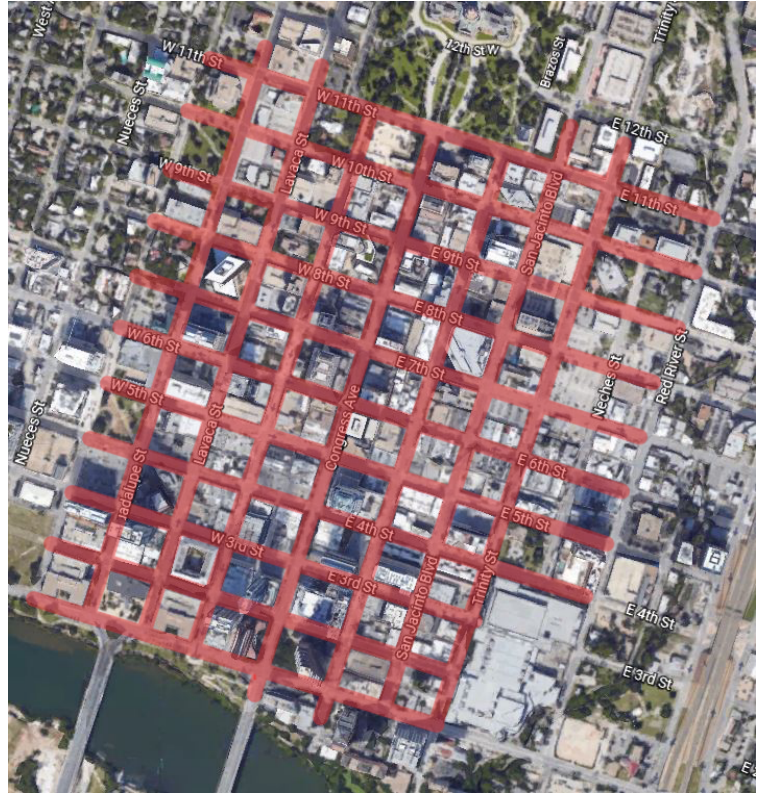

Figure 9: Austin downtown network (original source: Google Earth)

Table 2: Comparison of computational times (in seconds) for different simulation lengths in Austin downtown network using a time step of 1 second

| Simulation Horizon (s) | CTM | | LTM | | FLH | |
|---|---|---|---|---|---|---|
| | link model | node model | link model | node model | link model | node model |
| **200** | 1.411 | 0.485 | 0.359 | 0.520 | 0.354 | 0.527 |
| **500** | 3.582 | 1.159 | 0.938 | 1.368 | 0.856 | 1.309 |
| **1000** | 8.334 | 2.562 | 1.888 | 2.570 | 1.799 | 2.596 |

Table 3: Comparison of computational times (in seconds) for different simulation lengths in Austin downtown network using a time step of 5 seconds

| Simulation Horizon (s) | CTM | | LTM | | FLH | |
|---|---|---|---|---|---|---|
| | link model | node model | link model | node model | link model | node model |
| **200** | 0.066 | 0.110 | 0.074 | 0.106 | 0.079 | 0.113 |
| **500** | 0.164 | 0.260 | 0.177 | 0.266 | 0.169 | 0.257 |
| **1000** | 0.308 | 0.467 | 0.333 | 0.483 | 0.341 | 0.526 |

*5.3. Discussion*

An important difference among the discussed models is that, while at any point $(x, t)$, the solutions generated by the CTM and the FLH converge toward the solution of the LWR PDE, the solution generated by the LTM converges only in specific cases.

We illustrate this in Figure 10, where we present a scenario in which we consider a triangular fundamental diagram of parameters $k_c = 0.037 \ veh/m$, $k_{jam} = 0.1297 \ veh/m$, $u = 20 \ m/s$, and $w = 3.5 \ m/s$. We assume that the upstream half of a link is congested ($k_1 = 0.1297 \ veh/m$), while the downstream half of the link is free flow ($k_2 = 0.01 \ veh/m$). In this specific case, calculating $N$ at the point A (10,600) by using the LTM procedure (Newell's method) would yield: $N_A = \min\{N_D + (x_A - x_D) \cdot k_j; N_U\} = -51.88$. The correct solution, instead, would correspond to: $N_A = N_{U'} = -63.65$.

A common procedure to avoid expansion waves would be partitioning the network by splitting links wherever an expansion wave would occur, that is, wherever the density would decrease over space. However, dividing the link presents two issues.

- First, this requires a modification of the topology of the network, which becomes a function of the choice of the initial conditions. This is problematic in case of a-priori unknown initial conditions (for example in estimation, optimization or robust control problems where a large number of random initial conditions is drawn according to a certain distribution). This would increase the computational overhead before the actual computation process.

- Second, splitting the link in two or more links would increase the computational time to find the solution by a factor of two or more, since the demand and supplies at the boundaries would have to be derived for each split link.

The proposed FLH algorithm avoids these computational issues by imposing a minor computational penalty on the original LTM (three computations per time step instead of two computations per time step, while in the domain of influence of the initial conditions). Splitting a link in two for example would require four computations per time step. More than one split may be required, depending on the number of initial condition blocks and their configuration.

It would still be possible to apply Newell's method to derive the solution at the link boundaries as we did in previous examples. Because of that, the LTM would have a slightly increased computational time, comparable to the FLH.

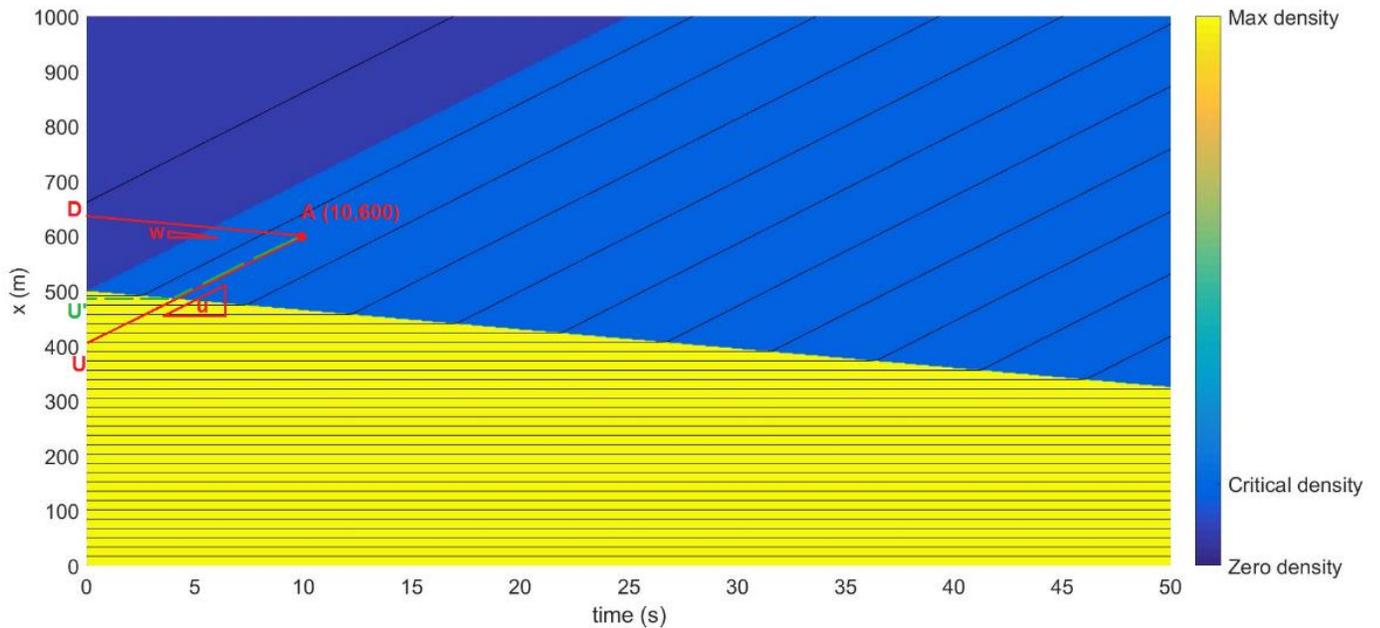

Figure 10: Derivation of the solution in point A by using Newell's method (solid lines) and correct approach (dashed line)

## 6. Conclusion

In this article, we introduce the Fast Lax Hopf (FLH) algorithm for solving the LWR-PDE on networks quickly and efficiently with any concave Fundamental Diagram. This algorithm presents the same accuracy of the original LH algorithm (it is exact in single link simulations), while improving the computational performance over the LH algorithm In particular, in the case of a triangular fundamental diagram, the performance of the FLH algorithm is comparable to the LTM (the LTM is the fastest solution method relying on time discretization to the best of authors' knowledge).

Like the CTM, but unlike the LTM, the proposed scheme converges everywhere (in the space-time domain) to the solution of the LWR model. The computational performance and accuracy of the FLH are analyzed in comparison with the original LH and two other known approaches: the LTM and the CTM. Clearly, depending on the typology of application and implementation details, each method can be the most suitable. The LTM is the fastest and guarantees a good level of accuracy, when initial conditions obey a specific structure (e.g. in the DNL) or for specific shapes of FDs. The FLH represents a good compromise in terms of accuracy, generality and computational speed. The CTM is perhaps the most "user-friendly" since it is easier to explain and implement (particularly for nontechnical audiences), and easier to customize (for example to solve multi-class or multi-commodity problems). Future work will be directed at applying the FLH in real-world traffic operations, and traffic control and estimation problems.

# Appendix I: Formulation of boundary and internal conditions for a general concave fundamental diagram

a) <u>Definition of initial, upstream, downstream and internal conditions</u>

The initial condition can be expressed as a piecewise linear function, with each linear piece defined by:

$$c_{ini_i}(x) = \begin{cases} -k_i x + b_i & : x_i \leq x \leq x_{i+1} \\ +\infty & : otherwise \end{cases} \tag{1}$$

With the above definition, the initial condition can be written as $c_{ini} = \min_i N_{ini_i}$

Similarly, the upstream boundary condition is assumed to be piecewise linear, with each piece defined by:

$$c_{up_j}(t) = \begin{cases} q_j t + d_j & : t_j \leq t \leq t_{j+1} \\ +\infty & : otherwise \end{cases} \tag{2}$$

With this definition, the upstream boundary condition can be written as $c_{up} = \min_j N_{up_j}$

The downstream boundary condition is also assumed to be a piecewise linear function, with each piece defined by:

$$c_{down_j}(t) = \begin{cases} p_j t + c_j & : t_j \leq t \leq t_{j+1} \\ +\infty & : otherwise \end{cases} \tag{3}$$

This enables us to define the downstream boundary condition function as $c_{down} = \min_j N_{down_j}$,

One of the major results of Mazaré et al. (2011) is that the solutions associated with each linear piece of the initial, upstream, downstream and internal boundary conditions can be computed analytically as follows:

b) <u>Solution to a linear initial condition</u>

Let $v_i \in \partial_+ Q(k_i)$   (physically, $v_i$ is the characteristic speed associated with $k_i$)

$$N_{c_{ini}}(x,t) = \begin{cases} -k_i \cdot x + Q(k_i) \cdot t + b_i & : x_i + tv_i \leq x \leq x_{i+1} + tv_i \\ k_i \cdot x_{i+1} + b_i + t \cdot R\left(\dfrac{x - x_{i+1}}{t}\right) & : x_{i+1} + tv_i \leq x \leq x_{i+1} + tv_f \\ k_i \cdot x_i + b_i + t \cdot R\left(\dfrac{x - x_i}{t}\right) & : x_i + tw \leq x \leq x_i + tv_i \end{cases} \tag{4}$$

c) <u>Solution to a linear upstream boundary condition</u>

For an upstream boundary condition $N_{up}$ defined as: $N_{up}{}^j(t) = q_j t + d_j$, let us define the associated density $\rho_j$ as $\rho_j = \inf_{\rho \in [0, k_{max}] \ s.t. \ Q(\rho) = q_j} \rho$ , and the associated characteristic speed $v_j \in \partial_+ Q(\rho_j)$. With this definition, the solution component can be expressed as:

$$N_{c_{up}{}^j}(x,t) = \begin{cases} d_j + \rho_j(x - x_0) + q_j \cdot t & : \ x_0 + v_j(t - t_{j+1}) \leq x \leq x_0 + v_j(t - t_j) \\ q_j t_j + d_j + (t - t_j)R\left(\dfrac{x - x_0}{t - t_j}\right) & : \ x_0 + v_j(t - t_j) \leq x \leq x_0 + v_f(t - t_j) \\ q_j t_{j+1} + d_j + (t - t_{j+1})R\left(\dfrac{x - x_0}{t - t_{j+1}}\right) & : \ x_0 \leq x \leq x_0 + v_j(t - t_{j+1}) \end{cases}$$ (6)

d) <u>Solution to a linear downstream boundary condition</u>

For a downstream boundary condition $N_{down}{}^j$, defined as $N_{down}{}^j(t) = p_j t + c_j$, let us define the associated density $\rho_j$ as $\rho_j = \sup_{\rho \in [0, k_{max}] \ s.t. \ Q(\rho) = q_j} \rho$ , and the associated characteristic speed $v_j \in \partial_+ Q(\rho_j)$. With this definition, the solution component can be expressed as:

$$N_{down}{}^j(x,t) = \begin{cases} t \cdot p_j + (x - x_n) \cdot \rho_j + c_j & : \ x_n + v_j(t - t_j) \leq x \leq v_j(t - t_{j+1}) \\ c_j + p_j t_j + (t - t_j)R\left(\dfrac{x - x_n}{t - t_j}\right) & : x_n + w(t - t_j) \leq x \leq x_n + v_j(t - t_j) \\ c_j + p_j t_{j+1} + (t - t_{j+1})R\left(\dfrac{x - x_n}{t - t_{j+1}}\right) & : x_n + v_j(t - t_{j+1}) \leq x \end{cases}$$ (7)

## Appendix II: Formulation of boundary and internal conditions for a triangular fundamental diagram

e) <u>Definition of initial, upstream, downstream and internal conditions</u>

The initial condition can be expressed as a piecewise linear function, with each linear piece defined by:

$$c_{ini_i}(x) = \begin{cases} -k_i x + b_i & : x_i \leq x \leq x_{i+1} \\ +\infty & : otherwise \end{cases} \tag{1}$$

With the above definition, the initial condition can be written as $c_{ini} = \min_i N_{ini_i}$

Similarly, the upstream boundary condition is assumed to be piecewise linear, with each piece defined by:

$$c_{up_j}(t) = \begin{cases} q_j t + d_j & : t_j \leq t \leq t_{j+1} \\ +\infty & : otherwise \end{cases} \tag{2}$$

With this definition, the upstream boundary condition can be written as $c_{up} = \min_j N_{up_j}$

The downstream boundary condition is also assumed to be a piecewise linear function, with each piece defined by:

$$c_{down_j}(t) = \begin{cases} p_j t + c_j & : t_j \leq t \leq t_{j+1} \\ +\infty & : otherwise \end{cases} \tag{3}$$

This enables us to define the downstream boundary condition function as $c_{down} = \min_j N_{down_j}$,

One of the major results of Mazaré et al. (2011) is that the solutions associated with each linear piece of the initial, upstream, downstream and internal boundary conditions can be computed analytically as follows:

f) <u>Solution to a linear initial condition</u>

If $0 \leq k_i \leq k_c$, the initial condition imposes a free-flow state.

$$N_{c_{ini}}(x,t) = \begin{cases} k_i\big(tv_f - x\big) + b_i & : x_i + tv_f \leq x \leq x_{i+1} + tv_f \\ k_c\big(tv_f - x\big) + b_i + x_i(k_c - k_i) & : x_i + tw \leq x \leq x_{i+1} + tv_f \end{cases} \tag{4}$$

else, if $k_c \leq k_i \leq k_j$, the initial condition imposes a congested state

$$N_{c_{ini_i}}(x,t)$$

$$= \begin{cases} k_i(tw - x) - tk_j w + b_i & : x_i + tw \leq x \leq x_{i+1} + tw \\ k_c(tw - x) - tk_j w + x_{i+1}(k_c - k_i) + b_i & : x_{i+1} + tw \leq x \leq x_{i+1} + tv_f \end{cases} \tag{5}$$

g) <u>Solution to a linear upstream boundary condition</u>

For an upstream boundary condition $N_{up}^{\,j}$ defined as: $N_{up_j}^{\,j}(t) = q_j t + d_j$ with $d_j = -q_j t + \sum_{l=0}^{j-1}(t_{l+1} - t_l)\,q_j^{\,l}$, the solution component can be expressed as:

$$N_{c_{up}}{}^{j}(x,t)$$

$$
= \begin{cases}
d_j + q_j\left(t - \dfrac{x - x_0}{v_f}\right) \; : \; x_0 + v_f(t - t_{j+1}) \le x \le x_0 + v_f(t - t_j) \\[2ex]
d_j + q_j t_{j+1} + k_c\left((t - t_{j+1})v_f - (x - x_0)\right) \; : \; x_0 \le x \le x_0 + v_f(t - t_{j+1})
\end{cases}
\tag{6}
$$

h)   <u>Solution to a linear downstream boundary condition</u>

For a downstream boundary condition $N_{down}{}^{j}$, defined as $N_{down}{}^{j}(t) = p_j t + b_j$ with $b_j = -p_j t + N_{ini}^{(n-1)}(x_n) + \sum_{l=0}^{j-1}(t_{l+1} - t_l)\, q_j{}^{l}$, the solution component can be expressed as:

$$N_{down}{}^{j}(x,t)$$

$$
= \begin{cases}
b_j + p_j t - \left(\dfrac{p_j}{w} + k_j\right)(x_n - x) \; : \; x_n + w(t - t_j) \le x \le x_n + w(t - t_{j+1}) \\[2ex]
b_j + p_j t_{j+1} + k_c\left((t - t_{j+1})v_f + x_n - x\right) \; : \; x_n + w(t - t_j) \le x \le x_n
\end{cases}
\tag{7}
$$

**Acknowledgements**

The authors would like to thank Chris Tampere for fruitful discussions on the Link Transmission Model. The authors would like to thank the Texas Department of Transportation for supporting this research under project 0-6838, Bringing Smart Transport to Texans: Ensuring the Benefits of a Connected and Autonomous Transport System in Texas. This research was supported by the National Science Foundation under award No 1636154.